\documentclass[11pt]{article}

\usepackage{amsmath,amssymb,amsthm}
\usepackage{hyperref}
\usepackage{geometry}
\usepackage{xcolor}
\usepackage{tikz}
\usetikzlibrary{backgrounds}

\newtheorem{theorem}{Theorem}[section]
\newtheorem{proposition}[theorem]{Proposition}
\newtheorem{lemma}[theorem]{Lemma}
\newtheorem{corollary}[theorem]{Corollary}
\newtheorem{conjecture}[theorem]{Conjecture}
\theoremstyle{definition}
\newtheorem{definition}[theorem]{Definition}
\newtheorem{remark}[theorem]{Remark}

\numberwithin{equation}{section}

\title{Asymptotics of linear normal lambda terms:\\
WKB structure, Airy turning point and a positive recurrence}
\author{Olivier Bodini}
\date{}

\begin{document}
\maketitle

\begin{abstract}
Let \(a_n\) be the number of closed linear \(\lambda\)-terms in
\(\beta\)-normal form with \(n\) variable occurrences. We prove
\[
\lim_{n\to\infty}\frac{a_n}{(n-1)!\,\lambda^n}=\frac{\kappa_\infty}{4}=\frac1{8\pi},
\qquad \lambda=\frac1{\ln2-\frac12},
\]
where \(\kappa_\infty=\lim_m K_m/(3^m\Gamma(m))\) is the large-order
constant of the Riccati--Airy sequence \(K_m\) attached to the turning point
of the WKB expansion. We show that the \(K_m\) are exactly rescaled
coefficients of the formal logarithmic derivative of the Airy asymptotic
series, and deduce \(\kappa_\infty=1/(2\pi)\) from the explicit Airy
coefficients. Hence
\[
a_n\sim\frac{(n-1)!}{8\pi}\Bigl(\frac1{\ln2-\frac12}\Bigr)^n .
\]
The proof combines four further ingredients:
(i) the WKB coefficients are polynomials with nonnegative coefficients in
\(v=(1-4w)^{-1/2}\) and satisfy an explicit positive recurrence; (ii) the
linearisation of the WKB equation has a combinatorial meaning
(\(\lambda\)-chain terms) and is a translation in the action coordinate
\(s\), whose singularity sits at the instanton action \(A_0=\ln2-\frac12\);
(iii) an exact identity shows that, after conjugation by the transported WKB
prefactor \(\mu=r/(1+r)^2\), the remaining non-linear terms are either
negative or of relative size \(O(m^{-2})\); (iv) under the conjugated linear
dynamics the depth profile is dominated by its top coefficient. A finite
bound on an explicit ratio sequence is obtained by exact computer
arithmetic (the program is given in an appendix); it also yields \(a_n\le1.14736\,g_n(1)\) for all \(n\ge2\), where
\(g_n\) is an explicit linear comparison class with
\(g_n(1)\sim\frac{5}{144}(n-1)!\lambda^n\).
\end{abstract}

\section{Introduction and main result}

A linear normal \(\lambda\)-term is either a free variable, an abstraction
\(\lambda x.M\) where \(x\) occurs exactly once in \(M\), or an application
\(M\,N\) where \(M\) is not an abstraction. Let \(L(z,u)\) count these terms,
\(z\) marking size and \(u\) free variables. Abstraction binds one free
variable (a derivative in \(u\)), and a term that is not an abstraction is a
spine \(x\,N_1\cdots N_p\) with a free head, so
\[
L=\frac{zu}{1-zL}+z\,\partial_uL .
\]
A term with \(n\) variable occurrences and \(k\) free variables has \(n-k\)
abstractions and \(n-1\) applications, hence size \(3n-k-1\). With
\(x=z^3\), \(t=u/z\) and \(L=F/z\) one gets
\begin{equation}\label{eq:core}
F(x,t)=\frac{xt}{1-F(x,t)}+\partial_tF(x,t),
\qquad
F=\sum_{n\ge1}f_n(t)\,x^n,\quad f_n(t)=\sum_{k=0}^na_{n,k}t^k,
\end{equation}
where \(a_{n,k}\) counts terms with \(n\) occurrences and \(k\) free
variables. We write \(a_n=a_{n,0}\) for closed terms (size \(3n-1\)):
\[
(a_1,\dots,a_7)=(1,\,3,\,26,\,367,\,7142,\,176766,\,5304356).
\]

\begin{theorem}\label{thm:intro}
Let \(\lambda=1/(\ln2-\frac12)\approx5.177399\) and let \(K_m\) be defined by
\(K_1=\frac12\) and
\(K_m=\sum_{i=1}^{m-1}K_iK_{m-i}+(3m-4)K_{m-1}\) \((m\ge2)\). Then
\(K_m/(3^m\Gamma(m))\to\kappa_\infty=\frac1{2\pi}\) and
\[
\lim_{n\to\infty}\frac{a_n}{(n-1)!\,\lambda^n}=\frac{\kappa_\infty}4=\frac1{8\pi},
\qquad\text{i.e.}\qquad
a_n\sim\frac{(n-1)!}{8\pi}\Bigl(\frac1{\ln2-\frac12}\Bigr)^n .
\]
\end{theorem}

The proof uses one computer-assisted input: a finite bound on an explicit
sequence \(U_m\), obtained by exact integer arithmetic
(Theorem~\ref{thm:Ubound}, program in Appendix~\ref{app:code}).

\paragraph{Related work.}
The analytic enumeration of \(\lambda\)-terms has been developed in
particular in a series of papers by Bodini and coauthors. Linear \(\lambda\)-terms
are the \(BCI(1)\)-terms of combinatory logic; Bodini, Gardy and Jacquot
\cite{BGJ13} counted \(BCI(1)\)- and \(BCK(1)\)-terms through a differential
equation for their generating function, obtained their asymptotic number and
designed random samplers. This was extended to \(BCI(p)\)-terms, with
differential equations also for \(BCK(p)\)-terms and general closed
\(\lambda\)-terms, in \cite{BGGJ13}, and to \(BCK(2)\)-terms in \cite{BG14}.
Terms with bounded unary height (nesting of abstractions) were studied in
\cite{BGG11,BGGG18}, statistical properties of random terms in de~Bruijn
notation in \cite{BBD19}, and limit laws of parameters of linear
\(\lambda\)-terms, via their bijection with rooted trivalent maps, in
\cite{BSZ25}. The present paper deals with linear terms in
\(\beta\)-normal form. The normal-form constraint replaces unrestricted
application by the spine term \(xt/(1-F)\) in \eqref{eq:core}; the resulting
equation has a WKB turning point, which governs both the growth rate
\(\lambda\) and the constant \(1/(8\pi)\).

\paragraph{Outline.}
Section~\ref{sec:wkb} sets up the WKB expansion in \(w=xt\) and
\(\varepsilon=x\), whose coefficients \(\varphi_m\) count terms with \(m\)
abstractions. Section~\ref{sec:K} proves \(\kappa_\infty=1/(2\pi)\) through an
exact reduction to the Airy series. Section~\ref{sec:chain} interprets the linearisation
(\(\lambda\)-chain terms) and introduces the action coordinate \(s\).
Section~\ref{sec:pos} proves positivity in \(v=1/r\) and an exact identity.
Section~\ref{sec:class} introduces the WKB-exact linear class.
Sections~\ref{sec:upper}--\ref{sec:main} prove the bound on \(U_m\) and the
main theorem. Section~\ref{sec:open} concludes and states three
conjectures.

\section{The WKB expansion}\label{sec:wkb}

\subsection{Structure in \(t\)}

Let \(H=1/(1-F)=\sum_{n\ge0}h_n(t)x^n\). Then
\(h_n=\sum_{k=1}^nf_kh_{n-k}\), \(h_0=1\), and \eqref{eq:core} reads
\(f_n'-f_n=-t\,h_{n-1}\).

\begin{lemma}\label{lem:t}
Each \(f_n\) is the polynomial \(f_n=\sum_{r\ge0}\partial_t^r(t\,h_{n-1})\),
with nonnegative coefficients, and \(a_{n,k}\le a_{n,0}/k!\). Hence
\(a_n\le f_n(t)\le e^t\,a_n\) for \(t\ge0\).
\end{lemma}

\begin{proof}
The combinatorial \(f_n\) is a polynomial, and the polynomial solution of
\(y'-y=-t\,h_{n-1}\) is unique (others differ by \(c\,e^t\)); it is
\(\sum_r\partial_t^r(th_{n-1})\), which has nonnegative coefficients by
induction. Comparing coefficients of \(t^k\),
\((k+1)a_{n,k+1}=a_{n,k}-[t^k](t\,h_{n-1})\le a_{n,k}\).
\end{proof}

\subsection{WKB coefficients}

Set \(w=xt\) and \(\varepsilon=x\). Then \(\partial_t=\varepsilon\partial_w\)
and \eqref{eq:core} becomes
\begin{equation}\label{eq:wkb}
\varepsilon\,\partial_wF=F-\frac{w}{1-F},
\qquad
F=\sum_{m\ge0}\varepsilon^m\varphi_m(w),\qquad
\varphi_m(w)=\sum_{k\ge0}a_{m+k,k}\,w^k .
\end{equation}
Since a term with \(n\) occurrences and \(k\) free variables has \(n-k\)
abstractions, \(\varphi_m\) is the generating function, by number of free
variables, of terms with exactly \(m\) abstractions, and
\(a_m=\varphi_m(0)\). Throughout, \(r=\sqrt{1-4w}\) and \('=d/dw\).

\begin{lemma}\label{lem:phi}
\(\varphi_0=\frac{1-r}2\), and for \(m\ge1\)
\begin{equation}\label{eq:star}
r\,\varphi_m=R_m:=\sum_{i=1}^{m-1}\varphi_i\varphi_{m-i}+\varphi'_{m-1}
-\sum_{i+j=m-1}\varphi_i\varphi'_j .
\end{equation}
Moreover \(\varphi_m=P_m(r)\,r^{-(3m-1)}\) with \(P_m\in\mathbb Q[r]\) and
\(P_m(0)=K_m>0\); in particular the pole at \(r=0\) has exact order
\(3m-1\).
\end{lemma}

\begin{proof}
Write \eqref{eq:wkb} as \(F(1-F)-w=\varepsilon(1-F)\partial_wF\). Order
\(0\) gives \(\varphi_0-\varphi_0^2=w\), \(\varphi_0(0)=0\). At order \(m\),
isolating \(\varphi_m\) and using \(1-2\varphi_0=r\) gives \eqref{eq:star}.
Since \(d/dw=-\frac2r\,d/dr\), a derivative raises the pole order by at
most \(2\). By induction (\(\varphi_1=\frac{1+r}{2r^2}\)), the terms of
\(R_m\) have pole order at most \(3m-2\), with exact order \(3m-2\) only for
\(\varphi_i\varphi_{m-i}\) (\(1\le i\le m-1\)), \(\varphi'_{m-1}\) and
\(-\varphi_0\varphi'_{m-1}\). Their leading coefficients are
\(K_iK_{m-i}\), \(2(3m-4)K_{m-1}\) and \(-(3m-4)K_{m-1}\), which gives the
recurrence of Theorem~\ref{thm:intro}, and \(K_m>0\) by induction.
\end{proof}

\section{The turning-point coefficients \(K_m\)}\label{sec:K}

\begin{theorem}\label{thm:kappa}
\(\displaystyle\kappa_m:=\frac{K_m}{3^m\Gamma(m)}=\frac1{2\pi}\bigl(1+O(1/m)\bigr)\).
In particular \(\kappa_\infty=\lim_m\kappa_m=1/(2\pi)\).
\end{theorem}

The proof relates \(K_m\) to the logarithm of the classical asymptotic
series of the Airy function \cite{Olver}, whose coefficients are explicit;
it is purely formal and uses no analytic property of Airy functions.

Let \(\Phi(z)=\sum_{m\ge1}K_mz^m\); the recurrence for \(K_m\) is equivalent
to
\begin{equation}\label{eq:riccatiPhi}
3z^2\Phi'=(1+z)\Phi-\Phi^2-\tfrac z2 .
\end{equation}
This is the local form of \eqref{eq:wkb} at the turning point: with
\(z=\varepsilon/r^3\), \(\eta=r^2/4\) and \(\tilde F=F-\frac12=r(\Phi(z)-\frac12)\), it
is the leading-order equation \(\varepsilon\,\partial_\eta\tilde F=2(\tilde F^2-\eta)\), which the
substitution \(\tilde F=-2^{-1/3}\varepsilon^{1/3}\psi'/\psi\) turns into the Airy
equation \(\psi''=\sigma\psi\).

\begin{lemma}[Exact Airy reduction]\label{lem:airyred}
Let \(Y(t)=\sum_{k\ge0}c_kt^k\), \(c_0=1\), be the formal power series such
that \(y(\sigma)=\sigma^{1/2}\,Y(\sigma^{-3/2})\) solves \(y'+y^2=\sigma\) (the
formal logarithmic derivative of the asymptotic expansion of the growing
Airy solution \(\mathrm{Bi}\), \(y=\sigma^{1/2}-\frac1{4\sigma}+\cdots\); the
normalisation \(c_0=1\) selects \(\mathrm{Bi}\), whereas \(c_0=-1\) would
correspond to the recessive solution \(\mathrm{Ai}\)). Then \(Y(t)=1-2\,\Phi(t/4)\), that
is, \(K_m=-\frac12\,4^m\,c_m\) for \(m\ge1\).
\end{lemma}

\begin{proof}
With \(t=\sigma^{-3/2}\), \(y'=\sigma^{-1/2}\bigl(\frac12Y-\frac32tY'\bigr)\), so
\(y'+y^2=\sigma\) is equivalent to
\begin{equation}\label{eq:riccatiY}
3t^2Y'=2Y^2-2+tY .
\end{equation}
Comparing coefficients of \(t^k\), \eqref{eq:riccatiY} has a unique formal
solution with \(c_0=1\)
(\(4c_k=(3k-4)c_{k-1}-2\sum_{i=1}^{k-1}c_ic_{k-i}\)). Substituting
\(Y(t)=1-2\Phi(t/4)\), \(t=4z\), turns \eqref{eq:riccatiY} into
\(-24z^2\Phi'=-8\Phi+8\Phi^2+4z-8z\Phi\), which is \eqref{eq:riccatiPhi};
the latter also has a unique solution with \(\Phi(0)=0\).
\end{proof}

\begin{lemma}[Airy coefficients]\label{lem:airycoef}
Let
\[
u_k=\frac{\Gamma(k+\frac16)\,\Gamma(k+\frac56)}{2\pi\,2^k\,k!}
=\prod_{j=1}^k\frac{(6j-5)(6j-1)}{72\,j},
\qquad S(\tau)=\sum_{k\ge0}u_k\tau^k,
\]
the coefficients of the asymptotic series of the Airy functions (\(u_0=1\));
in particular
\(\mathrm{Bi}(\sigma)\sim\pi^{-1/2}\sigma^{-1/4}e^{\zeta}S(1/\zeta)\) with
\(\zeta=\frac23\sigma^{3/2}\) \cite{Olver}.
Then, with \(\tau=\frac32t\) (that is, \(\tau=1/\zeta\)),
\[
Y(t)=1-\frac t4-\tau^2\,\frac{S'(\tau)}{S(\tau)} .
\]
Consequently, writing \(\log S(\tau)=\sum_{k\ge1}\ell^S_k\tau^k\),
\begin{equation}\label{eq:cn}
c_n=-\Bigl(\frac32\Bigr)^n(n-1)\,\ell^S_{n-1}\qquad(n\ge2).
\end{equation}
\end{lemma}

\begin{proof}
The two expressions of \(u_k\) agree since
\(\Gamma(\frac16)\Gamma(\frac56)=\pi/\sin\frac\pi6=2\pi\). The recurrence
\(72(k+1)u_{k+1}=(6k+1)(6k+5)u_k\) says that \(S\) is a formal solution of
\[
\mathcal L[S]:=36\tau^2S''+72(\tau-1)S'+5S=0 .
\]
Let \(Y_0(t)=1-\frac t4-\tau^2S'(\tau)/S(\tau)\). A direct computation gives
\[
\Bigl(3t^2Y_0'-2Y_0^2+2-tY_0\Bigr)\,S^2=-\frac{\tau^2}{18}\,S\,\mathcal L[S]=0,
\]
and \(S\) is invertible (\(u_0=1\)), so \(Y_0\) solves \eqref{eq:riccatiY}
with \(Y_0(0)=1\); hence \(Y_0=Y\). Finally
\(\tau^2S'/S=\tau^2(\log S)'=\sum_kk\,\ell^S_k\tau^{k+1}\), and the
coefficient of \(t^n\) is obtained with \(\tau^n=(\frac32)^nt^n\).
\end{proof}

\begin{lemma}\label{lem:logS}
\(\ell^S_n=u_n\bigl(1+O(1/n)\bigr)\).
\end{lemma}

\begin{proof}
For \(i\ge1\), \(u_i>0\) and \(u_i2^i/\Gamma(i)\to\frac1{2\pi}\); hence there are
constants \(0<c_1\le c_2\) with \(c_1\Gamma(i)2^{-i}\le u_i\le c_2\Gamma(i)2^{-i}\).
Since \(\sum_{i=1}^{N-1}\Gamma(i)\Gamma(N-i)=\Gamma(N-1)\sum_{j=0}^{N-2}\binom{N-2}j^{-1}\le4\,\Gamma(N-1)\)
for \(N\ge2\) (with \(p=N-2\), the last sum equals \(1\) if \(p=0\); if
\(p\ge1\), its two extreme terms equal \(1\) and its \(p-1\) other terms are at
most \(1/p\), so it is at most \(3\)), we get
\(\sum_{i=1}^{N-1}u_iu_{N-i}\le A\,u_{N-1}\) with \(A=2c_2^2/c_1\). Let
\(B=S-1\), so that the coefficients of \(B\) are nonnegative. By induction on
\(j\), \([\tau^n]B^j\le A^{j-1}u_{n-j+1}\) for \(n\ge j\ge1\). Since
\(\log S=\sum_{j\ge1}(-1)^{j+1}B^j/j\),
\[
|\ell^S_n-u_n|\le\sum_{j=2}^nA^{j-1}u_{n-j+1}
\le A\,u_{n-1}\sum_{j\ge2}A^{j-2}\frac{u_{n-j+1}}{u_{n-1}} .
\]
Since \(u_{k-1}/u_k=\frac{72k}{(6k-5)(6k-1)}\le\frac Ck\) with \(C=\frac{72}5\),
writing \(j=i+2\) we get
\[
\frac{u_{n-j+1}}{u_{n-1}}=\prod_{k=n-i}^{n-1}\frac{u_{k-1}}{u_k}
\le\frac{C^i}{(n-1)(n-2)\cdots(n-i)}\le\frac{C^i}{i!}
\qquad(0\le i\le n-2),
\]
so the last sum is at most \(\sum_{i\ge0}(AC)^i/i!=e^{AC}\), uniformly in
\(n\). Hence \(\ell^S_n=u_n+O(u_{n-1})=u_n(1+O(1/n))\), using again
\(u_{n-1}/u_n\le C/n\).
\end{proof}

\begin{proof}[Proof of Theorem~\ref{thm:kappa}]
By Lemmas~\ref{lem:airyred}, \ref{lem:airycoef} and~\ref{lem:logS}, for
\(m\ge2\),
\[
K_m=-\tfrac12\,4^mc_m=\tfrac12\,6^m\,(m-1)\,u_{m-1}\bigl(1+O(1/m)\bigr)
=\frac{6^m\,\Gamma(m-\frac56)\Gamma(m-\frac16)}{4\pi\,2^{m-1}\,(m-2)!}\bigl(1+O(1/m)\bigr).
\]
Since \(\Gamma(m-\frac56)\Gamma(m-\frac16)=\Gamma(m)\Gamma(m-1)(1+O(1/m))\), this
is \(\frac1{2\pi}3^m\Gamma(m)(1+O(1/m))\).
\end{proof}

\begin{remark}
Numerically \(m(1-2\pi\kappa_m)\to0.278\ldots\), in agreement with the
\(O(1/m)\) rate. The exponential factor \((\frac34)^m\) in
\(c_m\sim-\frac1\pi(\frac34)^m\Gamma(m)\) reflects the exponentially small
ratio \(\mathrm{Ai}/\mathrm{Bi}\approx e^{-\frac43\sigma^{3/2}}=e^{-4/(3t)}\);
the proof above avoids any Stokes or Borel argument because the logarithm of
\(S\) inherits its large-order behaviour directly from the explicit
coefficients \(u_k\).
\end{remark}

\section{Linearisation, \(\lambda\)-chain terms and the action coordinate}
\label{sec:chain}

Let \(c(w)=\frac2{1+r}=\sum_{n\ge0}C_nw^n\) be the Catalan series, so that
\(\varphi_0=w\,c\) counts \(\lambda\)-free terms, and put
\(\delta=F-\varphi_0\). Expanding \(\frac w{1-F}=\frac{wc}{1-c\delta}\)
according to the number \(k\) of arguments of the head spine that contain a
\(\lambda\) gives
\begin{equation}\label{eq:decomp}
\delta=M\Bigl[\varepsilon\,\partial_w(\varphi_0+\delta)+\mathcal N\Bigr],\qquad
\mathcal N=\sum_{k\ge2}w\,c^{k+1}\delta^k,\qquad
M=\frac1{1-wc^2}=\frac{1+r}{2r}.
\end{equation}
Here \(M=\sum_j(wc^2)^j\) counts \(\lambda\)-free contexts with one hole in
argument position, and \(\mathcal N\) counts \emph{branching spines}, i.e.\
spines with at least two arguments containing a \(\lambda\). Solving
\eqref{eq:decomp} for \(\delta'\) gives the exact equation
\(\varepsilon\delta'=\frac{2r}{1+r}\,\delta-\varepsilon\varphi_0'-\mathcal N\).
Its linearisation around \(\varphi_0\) is
\(\varepsilon\delta'=\frac{2r}{1+r}\,\delta-\varepsilon\varphi_0'\), whose
homogeneous part has coefficient \(1/M\); linearising thus amounts to
deleting \(\mathcal N\).

\begin{definition}
A linear normal term is a \emph{\(\lambda\)-chain term} if its abstraction
nodes are pairwise comparable for the ancestor relation.
\end{definition}

\begin{proposition}\label{prop:chain}
The series \(\varphi_0+\sum_{m\ge1}\varepsilon^m\varphi^{\rm lin}_m\),
\(\varphi^{\rm lin}_m=(M\partial_w)^m\varphi_0\), counts \(\lambda\)-chain
terms. In particular \(\ell_m:=\varphi^{\rm lin}_m(0)\le a_m\).
\end{proposition}

\begin{proof}
Two incomparable \(\lambda\)-nodes have an application node as lowest
common ancestor; it belongs to a spine whose head is a variable, so the two
\(\lambda\)-nodes lie in distinct arguments and the spine is branching.
Conversely a branching spine contains two incomparable \(\lambda\)-nodes.
The inequality follows from \eqref{eq:decomp}: \(\varphi_m=M(\varphi'_{m-1}+[\varepsilon^m]\mathcal N)\)
with \(M\) and \(\mathcal N\) having nonnegative coefficients, so
\(\varphi_m\gg\varphi^{\rm lin}_m\) by induction; evaluate at \(w=0\).
\end{proof}
\begin{figure}[ht]
\centering
\begin{tikzpicture}[
  level distance=11mm,
  level 1/.style={sibling distance=16mm},
  level 2/.style={sibling distance=16mm},
  level 3/.style={sibling distance=14mm},
  level 4/.style={sibling distance=12mm},
  every node/.style={font=\small},
  lam/.style={circle,draw=blue!70!black,fill=blue!10,inner sep=0.5pt,minimum size=7mm},
  app/.style={circle,draw,inner sep=0.5pt,minimum size=6mm},
  br/.style={circle,draw=red!70!black,fill=red!15,very thick,inner sep=0.5pt,minimum size=6mm},
  var/.style={inner sep=2pt},
  bind/.style={->,densely dashed,gray!70,thin}
]
\begin{scope}
\node[lam] (a1) {$\lambda x$}
  child {node[lam] (a2) {$\lambda y$}
    child {node[app] (a3) {$@$}
      child {node[var] (ax) {$x$}}
      child {node[lam] (a4) {$\lambda z$}
        child {node[app] (a5) {$@$}
          child {node[var] (az) {$z$}}
          child {node[var] (ay) {$y$}}}}}};
\begin{scope}[on background layer]
  \draw[line width=9pt,blue!12,line cap=round,line join=round]
    (a1.center)--(a2.center)--(a3.center)--(a4.center);
\end{scope}
\draw[bind] (ax) to[bend left=50] (a1);
\draw[bind] (ay) to[bend right=60] (a2);
\draw[bind] (az) to[bend left=50] (a4);
\node[align=center,below=3mm] at (0,-5.4) {\(\lambda x.\lambda y.\,x\,(\lambda z.\,z\,y)\)\\[2pt]
  \(\lambda\)-chain term};
\end{scope}
\begin{scope}[xshift=7cm]
\node[lam] (b1) {$\lambda x$}
  child {node[br] (b2) {$@$}
    child {node[app] (b3) {$@$}
      child {node[var] (bx) {$x$}}
      child {node[lam] (b4) {$\lambda y$}
        child {node[var] (by) {$y$}}}}
    child[sibling distance=22mm] {node[lam] (b5) {$\lambda z$}
      child {node[var] (bz) {$z$}}}};
\draw[bind] (bx) to[bend left=50] (b1);
\draw[bind] (by) to[bend right=50] (b4);
\draw[bind] (bz) to[bend right=50] (b5);
\node[align=center,below=3mm] at (0,-5.4) {\(\lambda x.\,x\,(\lambda y.y)\,(\lambda z.z)\)\\[2pt]
  not a \(\lambda\)-chain term};
\end{scope}
\end{tikzpicture}
\caption{Left: the abstraction nodes (blue) all lie on one root-to-leaf path
(shaded), so the term is a \(\lambda\)-chain term; every spine has at most
one argument containing a \(\lambda\). Right: \(\lambda y\) and \(\lambda z\)
are incomparable; their lowest common ancestor (red) belongs to the spine
\(x\,(\lambda y.y)\,(\lambda z.z)\), which has two arguments containing a
\(\lambda\). This is a \emph{branching spine}, counted by \(\mathcal N\) in
\eqref{eq:decomp} and deleted by the linearisation. It is the only closed linear normal
term with \(3\) variables that is not a \(\lambda\)-chain term. Dashed arrows
show the bindings (each \(\lambda\) binds exactly one occurrence); \(@\)
denotes application, with the function on the left.}
\label{fig:chain}
\end{figure}

Figure~\ref{fig:chain} illustrates the definition. Exhaustive enumeration
gives \((\ell_1,\dots,\ell_5)=(1,3,25,339,6361)\), against
\((a_1,\dots,a_5)=(1,3,26,367,7142)\).

\begin{lemma}\label{lem:asym}
The entire function \(h(r)=(1+r)\,e^{r^2/2-r}\) has \(r=0\) as its only
critical point, with critical value \(1\), and \(0\) as its only finite
asymptotic value.
\end{lemma}

\begin{proof}
\(h'(r)=r^2e^{r^2/2-r}\). Let \(\gamma(t)\to\infty\) (\(t\to1^-\)) with
\(h(\gamma(t))\to a\ne0\). Then
\(\log|h(\gamma)|=\ln|1+\gamma|+\mathrm{Re}(\gamma^2/2-\gamma)\) stays bounded,
so \(\mathrm{Re}(\gamma^2)=O(|\gamma|)\). Writing \(\gamma=\rho e^{i\theta}\),
this gives \(|\cos2\theta|\le C/\rho\). Since \(|\cos2\theta|\) is bounded
below away from the four directions \(\theta_0\in\{\pm\frac\pi4,\pm\frac{3\pi}4\}\),
and \(\gamma\) is continuous, for \(t\) close to \(1\) the path stays in a
sector \(S\) of small opening around one of them, at distance \(>2\) from
the origin. On \(S\), \(G(r)=\frac{r^2}2-r+\mathrm{Log}(1+r)\), with the
principal logarithm (\(1+r\) stays away from \((-\infty,0]\)), is holomorphic
and \(e^G=h\). Since \(h(\gamma)\to a\ne0\), \(G(\gamma(t))\) differs from a
fixed determination of \(\log a\) by \(o(1)+2\pi ik(t)\) with \(k\)
continuous and integer-valued, hence constant: \(G(\gamma(t))\) is bounded.
But \(|\mathrm{Im}\,G(\gamma)|\ge\frac{\rho^2}2|\sin2\theta|-\rho-\pi\to\infty\),
because \(|\sin2\theta|\to1\). This is a contradiction. Finally, \(0\) is an
asymptotic value: \(|h(iR)|=\sqrt{1+R^2}\,e^{-R^2/2}\to0\) as \(R\to+\infty\).
\end{proof}

\begin{theorem}\label{thm:action}
Let \(A(r)=\int_0^r\frac{u^2}{1+u}\,du=\frac{r^2}2-r+\ln(1+r)\),
\(A_0=A(1)=\ln2-\frac12\), and \(s(w)=A_0-A(r(w))\). Then
\(\frac{ds}{dw}=\frac{2r}{1+r}=\frac1M\), so that \(M\partial_w=\partial_s\),
and
\[
\Lambda(s):=\sum_{m\ge1}\ell_m\frac{s^m}{m!}=\frac{1-r(s)}2,
\qquad A(r(s))=A_0-s,\quad r(0)=1 .
\]
The branch \(r(s)\) extends analytically to the slit disk
\(\mathcal D=\{|s|<|A_0+2\pi i|\}\setminus[A_0,\infty)\), does not vanish
there, and near \(s=A_0\)
\begin{equation}\label{eq:puiseux}
r=v_0+\tfrac14v_0^2+O(v_0^3),\qquad v_0=\bigl(3(A_0-s)\bigr)^{1/3}.
\end{equation}
\end{theorem}

\begin{proof}
\(\frac{dA}{dr}=\frac{r^2}{1+r}\) and \(\frac{dr}{dw}=-\frac2r\) give
\(\frac{ds}{dw}=\frac{2r}{1+r}\). As \(s(w)=w+O(w^2)\) is invertible,
\(\varphi^{\rm lin}_m=\partial_s^m\varphi_0\), and \(w=0\) corresponds to
\(s=0\). Exponentiating, \(r(s)\) is the branch, with \(r(0)=1\), of the
inverse of \(h(r)=(1+r)e^{r^2/2-r}\) evaluated at \(e^{A_0-s}\) (it exists
near \(s=0\) since \(h'(1)\ne0\)).

The set of singular values of \(h\) (closure of the set of critical and
finite asymptotic values) is \(\{0,1\}\) by Lemma~\ref{lem:asym}. Hence, with
\(\Omega=\mathbb C\setminus h^{-1}(\{0,1\})\), the restriction
\(h:\Omega\to\mathbb C\setminus\{0,1\}\) is a covering map
\cite{Bergweiler}; it is also a local biholomorphism, because the only
critical point \(0\) lies in \(h^{-1}(1)\). Let \(E(s)=e^{A_0-s}\). It never
vanishes, and \(E(s)=1\) only for \(s\in A_0+2\pi i\mathbb Z\). None of these
points lies in \(\mathcal D\): \(A_0\) is on the slit, and
\(|A_0+2\pi ik|\ge|A_0+2\pi i|\) for \(k\ne0\). Thus \(E\) maps \(\mathcal D\)
into \(\mathbb C\setminus\{0,1\}\). The domain \(\mathcal D\) is star-shaped
with respect to \(0\) (a segment \([0,s]\) meets \([A_0,\infty)\) only if
\(s\) itself lies on this half-line), hence simply connected and locally
path-connected. Moreover \(h(1)=2e^{-1/2}=e^{A_0}=E(0)\), and \(1\in\Omega\)
since \(e^{A_0}\ne0,1\). By the lifting criterion for coverings, \(E\) has a
unique continuous lift \(\tilde r:\mathcal D\to\Omega\) with
\(h\circ\tilde r=E\) and \(\tilde r(0)=1\). Since \(h\) is a local
biholomorphism on \(\Omega\), \(\tilde r\) is holomorphic, and by uniqueness
of lifts on a small disk around \(0\) it coincides there with the branch
\(r(s)\). Hence \(r=\tilde r\) extends analytically to \(\mathcal D\), and it
does not vanish there because \(0\in h^{-1}(1)\) is not in \(\Omega\).
Finally, near
\(r=0\), \(A(r)=\frac{r^3}3\bigl(1-\frac34r+O(r^2)\bigr)\), so
\(v_0=(3A(r))^{1/3}=r\bigl(1-\frac14r+O(r^2)\bigr)\) is a local
biholomorphism; along \([0,A_0)\) the function \(r\) is real and positive,
which selects the real cube root. Inverting gives \eqref{eq:puiseux}.
\end{proof}

\begin{corollary}\label{cor:chain}
\(\ell_n\sim c_\ell\,(n-1)!\,\lambda^n\,n^{-1/3}\) with
\(c_\ell=(3A_0)^{1/3}/(2|\Gamma(-\frac13)|)\approx0.10261\). Thus the
instanton action \(A_0\) is the radius of the Borel transform of the
linearised class, whose exponent is \(-\frac13\).
\end{corollary}

\begin{proof}
By \eqref{eq:puiseux}, \(\Lambda=\frac12-\frac12v_0+O(v_0^2)\); apply the
transfer theorems of singularity analysis \cite{FS} on the \(\Delta\)-domain
provided by Theorem~\ref{thm:action}.
\end{proof}

\section{Positivity in \(v=1/r\) and an exact identity}\label{sec:pos}

Put \(v=1/r=(1-4w)^{-1/2}\). Then \(\frac d{dw}=2v^3\frac d{dv}\),
\(M=\frac{1+v}2\), \(\varphi_0=\frac12(1-v^{-1})\), \(\varphi_0'=v\), and
\[
T:=M\partial_w=v^3(1+v)\,\partial_v=\partial_s .
\]
Since \(v^a\) has nonnegative coefficients in \(w\) for \(a>0\), we write
\(\psi_1\ll\psi_2\) for coefficientwise inequalities between polynomials in \(v\);
they imply the same inequalities in \(w\). For \(m\ge1\) we write
\(\varphi_m=\sum_d\varphi_{m,d}\,v^{3m-1-d}\) (\(d\) is the \emph{depth}
below the turning point), so that \(a_m=\varphi_m(v{=}1)=\sum_d\varphi_{m,d}\), and
\[
q_{m,d}=\sum_{i=1}^{m-1}\sum_{a+b=d}\varphi_{i,a}\varphi_{m-i,b},\qquad
q^{\rm mid}_{m,d}=\sum_{i=2}^{m-2}\sum_{a+b=d}\varphi_{i,a}\varphi_{m-i,b}
\]
for the depth profiles of \(v\sum_{i=1}^{m-1}\varphi_i\varphi_{m-i}\) and of
its part without the two edge terms \(i=1,m-1\).

\begin{lemma}\label{lem:pos}
For \(m\ge1\), \(\varphi_{m,d}=0\) unless \(0\le d\le m\), and
\(\varphi_{m,d}\ge0\). More precisely, \(\varphi_{1,\cdot}=(\frac12,\frac12)\) and
for \(m\ge2\)
\begin{equation}\label{eq:posrec}
\varphi_{m,d}=(3m-4-d)\bigl(\varphi_{m-1,d}+\varphi_{m-1,d-1}-q_{m-1,d-1}\bigr)+q_{m,d},
\end{equation}
where \(\chi_{m,d}:=\varphi_{m,d}-q_{m,d}\) satisfies
\(\chi_{m,d}=(3m-4-d)(\varphi_{m-1,d}+\chi_{m-1,d-1})\ge0\).
\end{lemma}

\begin{proof}
By \eqref{eq:star}, \(\varphi_m=v\,R_m\). In \(R_m\), the terms containing
\(\varphi_0\) form the linear part
\((\frac12+\frac1{2v})\varphi'_{m-1}-v\varphi_{m-1}\); after multiplication by
\(v\), it sends \(v^a\) to \(a\,v^{a+3}+(a-1)\,v^{a+2}\), that is
\((3m-4-d)(\varphi_{m-1,d}+\varphi_{m-1,d-1})\) at depth \(d\). The term
\(v\sum\varphi_i\varphi_{m-i}\) gives \(q_{m,d}\). In
\(-v\sum_{i,j\ge1,\,i+j=m-1}\varphi_i\varphi'_j\) the weights
\(2(3j-1-b)\) symmetrise to \((3j-1-b)+(3i-1-a)=3m-4-d\), giving
\(-(3m-4-d)q_{m-1,d-1}\). Positivity follows by induction since
\(3m-4-d\ge0\) for \(d\le m\).
\end{proof}

\begin{lemma}\label{lem:exactid}
For \(m\ge4\) and \(0\le d\le m\),
\begin{equation}\label{eq:exactid}
\varphi_{m,d}=(3m-3-d)\,\varphi_{m-1,d}+(3m-4-d)\,\varphi_{m-1,d-1}-\mathcal E_{m,d}+q^{\rm mid}_{m,d},
\end{equation}
\[
\mathcal E_{m,d}=\varphi_{m-2,d-1}+\varphi_{m-2,d-2}+(3m-6-d)\,q_{m-2,d-2}
+(3m-5-d)\,q^{\rm mid}_{m-1,d-1}\ \ge0 .
\]
\end{lemma}

\begin{proof}
As \(\varphi_{1,\cdot}=(\frac12,\frac12)\), the two edge terms of \(q_{n,e}\)
(distinct for \(n\ge3\)) give
\(q_{n,e}=\varphi_{n-1,e}+\varphi_{n-1,e-1}+q^{\rm mid}_{n,e}\). Inserting
this into \eqref{eq:posrec} gives, for \(m\ge3\),
\begin{equation}\label{eq:E1}
\varphi_{m,d}=(3m-3-d)(\varphi_{m-1,d}+\varphi_{m-1,d-1})-(3m-4-d)\,q_{m-1,d-1}+q^{\rm mid}_{m,d}.
\end{equation}
For \(m\ge4\), expand \(q_{m-1,d-1}\) in the same way and use \eqref{eq:E1}
at rank \(m-1\) to write
\((3m-4-d)(\varphi_{m-2,d-1}+\varphi_{m-2,d-2})=\varphi_{m-1,d-1}+\varphi_{m-2,d-1}+\varphi_{m-2,d-2}+(3m-6-d)\,q_{m-2,d-2}-q^{\rm mid}_{m-1,d-1}\).
Collecting terms gives \eqref{eq:exactid}. The coefficients
\(3m-6-d\ge2m-6\) and \(3m-5-d\ge2m-5\) are positive.
\end{proof}

At \(d=0\), \eqref{eq:exactid} reads
\(K_m=(3m-3)K_{m-1}+\sum_{i=2}^{m-2}K_iK_{m-i}\): it is the cancellation in
the proof of Theorem~\ref{thm:kappa}, now holding at every depth. The edge
terms correspond to branchings whose side branch carries a single
abstraction.

\section{The WKB-exact linear class}\label{sec:class}

\begin{proposition}\label{prop:class}
Let \(T''=v^3(1+v)\partial_v+v^3-v^2\). For a polynomial \(\psi\) of top
degree \(3(m-1)-1\), \([T''\psi]_d=(3m-3-d)\psi_d+(3m-4-d)\psi_{d-1}\), so
\eqref{eq:exactid} reads
\(\varphi_m=T''\varphi_{m-1}-\mathcal E_m+v\sum_{i=2}^{m-2}\varphi_i\varphi_{m-i}\).
Moreover
\[
T''=\mu\,\partial_s\,\mu^{-1},\qquad \mu=\frac{r}{(1+r)^2}=\frac{v}{(1+v)^2},
\]
where \(\mu\) is the WKB prefactor transported along the action. Let
\(g_2=\varphi_2=\frac54v^5+\frac32v^4+\frac14v^3\) and
\(g_m=T''g_{m-1}\). Then \(g_{m,0}=\frac54\,3^{m-2}(m-1)!\),
\(g_{m,d}>0\) for \(0\le d\le m\), and
\[
g_m(1)=\tfrac14\,(m-2)!\,[s^{m-2}]\frac{\varphi_2}{\mu}(s)
=\frac5{144}\,(m-1)!\,\lambda^m\Bigl(1+c_f\,m^{-1/3}+O(m^{-2/3})\Bigr),
\]
with
\(c_f=\frac{17}{8}\cdot\frac{36}5\cdot3^{-5/3}A_0^{1/3}/\Gamma(\frac53)\approx1.570\).
\end{proposition}

\begin{proof}
\(T''v^a=(a+1)v^{a+3}+(a-1)v^{a+2}\) gives the depth formula. Since
\(\mu^{-1}=(1+v)^2/v\), one has
\(T\mu^{-1}/\mu^{-1}=v^3(1+v)\bigl(\frac2{1+v}-\frac1v\bigr)=v^3-v^2\), hence
\(T(\mu^{-1}\psi)=\mu^{-1}\,T''\psi\), and \(T=\partial_s\). Thus
\(\mu^{-1}g_m=\partial_s^{m-2}(\mu^{-1}\varphi_2)\). At \(v=1\) (that is, \(s=0\))
\(\mu^{-1}=4\), which gives the coefficient formula with
\(\varphi_2/\mu=(1+v)^2(\frac54v^4+\frac32v^3+\frac14v^2)\). By
Theorem~\ref{thm:action}, \(v=1/r(s)\) is analytic in \(\mathcal D\) and
\(v=v_0^{-1}(1-\frac14v_0+O(v_0^2))\), so
\(\varphi_2/\mu=\frac54v_0^{-6}+\frac{17}8v_0^{-5}+O(v_0^{-4})\). Transfer
\cite{FS} with \([s^n](A_0-s)^{-2}=(n+1)A_0^{-n-2}\) gives the expansion.
Positivity follows by induction on \(m\ge2\), starting from the positive
coefficients \(\frac54,\frac32,\frac14\) of \(g_2\): by the depth formula,
\(g_{m,d}=(3m-3-d)g_{m-1,d}+(3m-4-d)g_{m-1,d-1}\) with nonnegative
multipliers for \(d\le m\); for \(0\le d\le m-1\), \(g_{m,d}\ge(3m-3-d)g_{m-1,d}>0\)
since \(3m-3-d\ge2m-2>0\), and for \(d=m\), \(g_{m-1,m}=0\) and
\(g_{m,m}=(2m-4)g_{m-1,m-1}>0\) since \(m\ge3\).
\end{proof}

The \(\lambda\)-chain class is generated by \(T=\partial_s\) and has
exponent \(-\frac13\); the conjugated class \(g\) has exponent \(0\)
because \(\varphi_2/\mu\) has a pole at the turning point. The heuristic
WKB transport is thus an exact conjugation, and \(\mu=\frac14\) at \(s=0\).

\begin{lemma}[Top dominance]\label{lem:topdom}
For \(a\ge1\) let \(J_n(a)=[T''^{\,n}v^a](1)\). Then
\[
J_n(a)=\tfrac14\,n!\,[s^n]\bigl((1+v)^2v^{a-1}\bigr)\ \sim\
\frac{(3A_0)^{-(a+1)/3}}{4\,\Gamma(\frac{a+1}3)}\;n!\,A_0^{-n}\,n^{\frac{a+1}3-1}.
\]
Hence, for fixed \(M\ge2\) and \(\psi=\sum_{d=0}^M\psi_dv^{3M-1-d}\) with
\(\psi_d\ge0\),
\[
\lim_{n\to\infty}\frac{[T''^{\,n}\psi](1)}{g_{M+n}(1)}=\frac{\psi_0}{g_{M,0}} .
\]
\end{lemma}

\begin{proof}
The dominant singular term of \((1+v)^2v^{a-1}\) at \(s=A_0\) is
\(v^{a+1}\sim(3(A_0-s))^{-(a+1)/3}\); apply transfer. Thus
\(J_n(a)/J_n(a')\to0\) for \(a<a'\), and both
\([T''^n\psi](1)\) and \(g_{M+n}(1)=[T''^ng_M](1)\) are finite sums
\(\sum_dc_dJ_n(3M-1-d)\) with \(c_d\ge0\). The denominator is
\(\asymp J_n(3M-1)\) since \(g_{M,0}>0\), and the numerator equals
\(\psi_0J_n(3M-1)+o(J_n(3M-1))\); this also covers \(\psi_0=0\), where the
limit is \(0\).
\end{proof}

\begin{lemma}\label{lem:profiles}
Let \(\hat E_m(x)=\sum_dg_{m,d}x^d/g_{m,0}\). Then
\(\hat E_2=1+\frac65x+\frac15x^2\ll(1+x)^2\), and for \(m\ge3\)
\(\hat E_m\ll(1+x)\hat E_{m-1}\ll\gamma_m\hat E_m\) with
\(\gamma_m=\frac{3m-3}{2m-4}=\frac32+\frac3{2m-4}\). Consequently, for
\(i,j\ge2\) and \(m=i+j\),
\(\hat E_i\hat E_j\ll C_{i,j}\hat E_m\) with
\(C_{i,j}=C_{j,i}=\prod_{k=\max(i,j)+1}^m\gamma_k\) (a product of
\(\min(i,j)\) factors).
\end{lemma}

\begin{proof}
\(\hat g_{m,d}=(1-\frac d{3m-3})\hat g_{m-1,d}+(1-\frac{d+1}{3m-3})\hat g_{m-1,d-1}\);
both multipliers lie in \([0,1]\) and are at least \(\frac{2m-4}{3m-3}\) on
the supports (\(d\le m-1\), resp.\ \(d\le m\)). Then
\(\hat E_i\ll(1+x)^{i-2}\hat E_2\ll(1+x)^i\). By symmetry we may assume
\(i\le j\); then \(\hat E_i\hat E_j\ll(1+x)^i\hat E_j\ll C_{i,j}\hat E_m\) by
\(i\) applications of the second inequality, at ranks \(j+1,\dots,m\).
\end{proof}

\section{A bound on the depth ratio}\label{sec:upper}

Let \(U_m=\max_{0\le d\le m}\varphi_{m,d}/g_{m,d}\) for \(m\ge2\), and
\[
k_{i,m}=\frac{g_{i,0}g_{m-i,0}}{g_{m,0}}=\frac{5/36}{(m-1)\binom{m-2}{i-1}},
\qquad
\varepsilon_m=\sum_{i=2}^{m-2}C_{i,m-i}\,k_{i,m}.
\]

\begin{proposition}\label{prop:Urec}
\(a_m\le U_m\,g_m(1)\). For \(m\ge4\),
\(q^{\rm mid}_m\ll\bigl(\sum_{i=2}^{m-2}U_iU_{m-i}C_{i,m-i}k_{i,m}\bigr)g_m\) and
\(U_m\le U_{m-1}+\sum_{i=2}^{m-2}U_iU_{m-i}\,C_{i,m-i}\,k_{i,m}\).
\end{proposition}

\begin{proof}
First, \(a_m=\sum_d\varphi_{m,d}\le U_m\sum_dg_{m,d}=U_m\,g_m(1)\). Let
\(m\ge4\). For \(j\ge2\) and every \(a\), the definitions of \(U_j\) and of
\(\hat E_j\) give \(\varphi_{j,a}\le U_jg_{j,a}=U_j\,g_{j,0}\,[x^a]\hat E_j(x)\).
For \(2\le i\le m-2\) both indices \(i\) and \(m-i\) are at least \(2\), so
\[
q^{\rm mid}_{m,d}=\sum_{i=2}^{m-2}\sum_{a+b=d}\varphi_{i,a}\varphi_{m-i,b}
\le\sum_{i=2}^{m-2}U_iU_{m-i}\,g_{i,0}\,g_{m-i,0}\,[x^d]\bigl(\hat E_i\hat E_{m-i}\bigr).
\]
By Lemma~\ref{lem:profiles},
\([x^d](\hat E_i\hat E_{m-i})\le C_{i,m-i}\,[x^d]\hat E_m=C_{i,m-i}\,g_{m,d}/g_{m,0}\),
and \(g_{i,0}g_{m-i,0}/g_{m,0}=k_{i,m}\) by Proposition~\ref{prop:class}. Hence
\[
q^{\rm mid}_{m,d}\le\Bigl(\sum_{i=2}^{m-2}U_iU_{m-i}\,C_{i,m-i}\,k_{i,m}\Bigr)g_{m,d}
\qquad(0\le d\le m),
\]
which is the second claim. Finally, dropping \(\mathcal E_m\ge0\) in
\eqref{eq:exactid} gives \(\varphi_m\ll T''\varphi_{m-1}+q^{\rm mid}_m\). The
depth multipliers \(3m-3-d\) and \(3m-4-d\) of \(T''\) are nonnegative for
\(d\le m\), so \(T''\varphi_{m-1}\ll U_{m-1}T''g_{m-1}=U_{m-1}g_m\). Dividing
by \(g_{m,d}>0\) and taking the maximum over \(d\) gives the bound on \(U_m\).
\end{proof}

\begin{theorem}[computer-assisted]\label{thm:Ubound}
\(U_m\le U^*:=1.14736\) for all \(m\ge2\). Consequently
\(a_n\le1.14736\,g_n(1)\) for all \(n\ge2\).
\end{theorem}

\begin{proof}
\emph{Exact part.} For \(m\le700\), the rational number \(U_m\) is computed
exactly with the program of Appendix~\ref{app:code}. The maximum is attained at \(d=0\) for
\(5\le m\le700\), and \(\bar U:=\max_{m\le700}U_m=U_{700}=1.1454596\ldots\)

\emph{Tail.} Let \(m\ge701\) and \(T_n=C_{n,m-n}k_{n,m}\), so that
\(\varepsilon_m\le2\sum_{2\le n\le m/2}T_n\) (the summand of \(\varepsilon_m\) is
invariant under \(i\mapsto m-i\)). Every numerical inequality below is an
inequality between explicit rational numbers; all of them are checked in
exact rational arithmetic by the program of Appendix~\ref{app:code}.

\emph{Region \(2\le n\le m/5\).} For \(n+1\le m/5\) one has
\(m-n\ge\frac45m+1>561\) and \(\frac n{m-1-n}\le\frac{m/5-1}{4m/5}<\frac14\);
since \(\gamma_k\) decreases in \(k\),
\(T_{n+1}/T_n=\gamma_{m-n}\frac n{m-1-n}\le\gamma_{561}/4\le0.37568\).
Moreover, as \(m-1\ge700\),
\(T_2=\gamma_{m-1}\gamma_m\frac{5/36}{(m-1)(m-2)}\le\frac5{36}\gamma_{700}\gamma_{701}\frac1{(m-1)(m-2)}\le\frac{0.31340}{(m-1)(m-2)}\).
Hence \(2\sum_{2\le n\le m/5}T_n\le\frac{2\cdot0.31340}{1-0.37568}\cdot\frac1{(m-1)(m-2)}
\le\frac{1.0040}{(m-1)(m-2)}\).

\emph{Region \(m/5<n\le m/2\).} Here \(C_{n,m-n}=\prod_{k=m-n+1}^m\gamma_k\)
with \(m-n+1>351\), so \(C_{n,m-n}\le\gamma_{351}^{\,n}\le1.51^{m/2}\).
Since \(\lfloor m/5\rfloor\le n-1\le(m-2)/2\), unimodality of binomial
coefficients gives
\(\binom{m-2}{n-1}\ge\binom{m-2}{\lfloor m/5\rfloor}\ge\bigl(\frac{m-2}{\lfloor m/5\rfloor}\bigr)^{\lfloor m/5\rfloor}\ge4.98^{\,m/5-1}\),
because \(\frac{m-2}{\lfloor m/5\rfloor}\ge5-\frac{10}m\ge5-\frac{10}{701}\ge4.98\).
So \(k_{n,m}\le\frac5{36}\,\frac{4.98^{-(m/5-1)}}{m-1}\). There are at most
\(\frac{3m}{10}+1\le\frac m2\) such \(n\); with the factor \(2\) in the bound on
\(\varepsilon_m\), their contribution is
\[
R_m\le\frac5{36}\cdot4.98\cdot\frac m{m-1}\Bigl(\frac{1.51^{5/2}}{4.98}\Bigr)^{m/5}\le0.70\,q^m,
\qquad q=0.8914,
\]
using \(\frac m{m-1}\le\frac{701}{700}\),
\(\frac5{36}\cdot4.98\cdot\frac{701}{700}\le0.70\) and
\(1.51^{5}\le(4.98\,q^5)^2\), i.e.\ \(\frac{1.51^{5/2}}{4.98}\le q^5\).

\emph{Summation.} Hence
\(\varepsilon_m\le\frac{1.0040}{(m-1)(m-2)}+0.70\,q^m\), and, using
\(\sum_{m\ge701}\frac1{(m-1)(m-2)}=\frac1{699}\) and
\(\sum_{m\ge701}0.70\,q^m=\frac{0.70\,q^{701}}{1-q}\le10^{-34}\),
\[
\eta:=\sum_{m>700}\varepsilon_m\le\frac{1.0040}{699}+10^{-34}\le0.001437 .
\]

\emph{Bootstrap.} Put \(U^*=1.14736\) and \(\eta=0.001437\). The exact value
of \(\bar U\) satisfies
\[
\bar U+U^{*2}\eta\le U^*
\]
(numerically \(\bar U+U^{*2}\eta=1.1473513\ldots\)). We prove \(U_m\le U^*\) for
all \(m\ge2\) by strong induction. For \(m\le700\), \(U_m\le\bar U\le U^*\).
Let \(m\ge701\) and assume \(U_k\le U^*\) for \(2\le k<m\). For
\(701\le k\le m\), the bound of Proposition~\ref{prop:Urec} on \(U_k\) involves
\(U_{k-1}\) and \(U_i,U_{k-i}\) with \(2\le i\le k-2\), all of index \(<k\); hence
\(U_k\le U_{k-1}+U^{*2}\varepsilon_k\). Summing over \(701\le k\le m\),
\(U_m\le U_{700}+U^{*2}\sum_{k=701}^m\varepsilon_k\le\bar U+U^{*2}\eta\le U^*\).
\end{proof}

\section{Exact asymptotics}\label{sec:main}

Let \(U=\sup_mU_m\le U^*\). The key point is to compare the remainder
\(\mathcal E_k\) with the reference profile \(g_k\) rather than with \(\varphi_k\):
\(\varphi_{k,d}\) is much smaller than \(g_{k,d}\) in the deepest layers
\(d\approx k\), but \(\mathcal E_k\) is uniformly small with respect to \(g_k\).

\begin{lemma}\label{lem:Rbound}
For \(k\ge4\), \(q^{\rm mid}_k\ll U^2\varepsilon_kg_k\). For \(k\ge9\),
\(\mathcal E_k\ll\bar\rho_kg_k\), where
\[
\bar\rho_k=U\Bigl[\tfrac1{(2k-3)(2k-4)}+\tfrac1{(2k-4)(2k-6)}+\tfrac1{(2k-6)(2k-3)}+\tfrac1{(2k-8)(2k-4)}\Bigr]
+U^2\Bigl[\tfrac{\varepsilon_{k-2}}{2k-4}+\varepsilon_{k-1}\Bigr].
\]
In particular \(\sum_k\varepsilon_k<\infty\) and \(\sum_k\bar\rho_k<\infty\).
\end{lemma}

\begin{proof}
The first claim is Proposition~\ref{prop:Urec}. A
horizontal step gives \(g_{k,d}\ge(3k-3-d)g_{k-1,d}\), a diagonal step
\(g_{k,d}\ge(3k-4-d)g_{k-1,d-1}\). Using \(\varphi_j\ll Ug_j\) (\(j\ge2\)) and the
edge decomposition of \(q_{k-2}\):
\begin{itemize}
\item \(\varphi_{k-2,d-1}\) (nonzero only if \(d\le k-1\)):
\(g_{k,d}\ge(3k-4-d)(3k-5-d)g_{k-2,d-1}\ge(2k-3)(2k-4)g_{k-2,d-1}\).
\item \(\varphi_{k-2,d-2}\):
\(g_{k,d}\ge(3k-4-d)(3k-6-d)g_{k-2,d-2}\ge(2k-4)(2k-6)g_{k-2,d-2}\).
\item \((3k-6-d)q_{k-2,d-2}\): the edge term \(\varphi_{k-3,d-2}\) contributes
at most \(\frac U{(3k-7-d)(3k-4-d)}\le\frac U{(2k-6)(2k-3)}\), the edge term
\(\varphi_{k-3,d-3}\) at most
\(\frac U{(3k-8-d)(3k-4-d)}\le\frac U{(2k-8)(2k-4)}\), and the middle part at
most \(\frac{U^2\varepsilon_{k-2}}{3k-4-d}\le\frac{U^2\varepsilon_{k-2}}{2k-4}\).
\item \((3k-5-d)q^{\rm mid}_{k-1,d-1}\le(3k-5-d)U^2\varepsilon_{k-1}g_{k-1,d-1}\le U^2\varepsilon_{k-1}g_{k,d}\).
\end{itemize}
Summability follows from the tail bound in the proof of
Theorem~\ref{thm:Ubound}.
\end{proof}

\begin{theorem}\label{thm:main}
\[
\lim_{n\to\infty}\frac{a_n}{g_n(1)}=\frac{36}5\,\kappa_\infty,
\qquad\text{hence}\qquad
\lim_{n\to\infty}\frac{a_n}{(n-1)!\,\lambda^n}=\frac{\kappa_\infty}4=\mu(0)\,\kappa_\infty .
\]
\end{theorem}

\begin{proof}
By \eqref{eq:exactid},
\(T''\varphi_{k-1}-\mathcal E_k\ll\varphi_k\ll T''\varphi_{k-1}+q^{\rm mid}_k\) for
\(k\ge4\). Since \(T''\) has nonnegative coefficients, induction from a
fixed \(M\ge8\), Lemma~\ref{lem:Rbound} and \(T''^{\,m-k}g_k=g_m\) give,
for \(m>M\),
\[
T''^{\,m-M}\varphi_M-P_M\,g_m\ \ll\ \varphi_m\ \ll\ T''^{\,m-M}\varphi_M+S_M\,g_m,
\qquad P_M=\sum_{k>M}\bar\rho_k,\quad S_M=U^2\sum_{k>M}\varepsilon_k .
\]
Evaluate at \(v=1\), divide by \(g_m(1)\), and let \(m\to\infty\) with
Lemma~\ref{lem:topdom} (\(\psi=\varphi_M\)):
\[
t_M-P_M\le\liminf_m\frac{a_m}{g_m(1)}\le\limsup_m\frac{a_m}{g_m(1)}\le t_M+S_M,
\qquad t_M=\frac{K_M}{g_{M,0}}=\frac{36}5\,\kappa_M .
\]
Let \(M\to\infty\): \(P_M,S_M\to0\) and \(t_M\to\frac{36}5\kappa_\infty\).
Conclude with Proposition~\ref{prop:class}, since
\(\frac5{144}\cdot\frac{36}5=\frac14\).
\end{proof}

\begin{corollary}\label{cor:main}
\[
a_n\sim\frac{(n-1)!}{8\pi}\Bigl(\frac1{\ln2-\frac12}\Bigr)^n .
\]
For every \(t\ge0\), the Borel transform
\(\mathcal Y(s,t)=\sum_nf_n(t)s^n/n!\) has radius of convergence \(A_0\)
(Lemma~\ref{lem:t}). Moreover
\[
\mathcal Y(s,0)\sim-\frac1{8\pi}\log(1-\lambda s)\qquad(s\to A_0^-).
\]
\end{corollary}

\begin{remark}[Speed of convergence]\label{rem:speed}
The convergence is slow:
\begin{center}
\begin{tabular}{r|cc}
\(n\) & \(a_n/g_n(1)\) & \(a_n/((n-1)!\lambda^n)\)\\ \hline
100 & 1.08074 & 0.05164\\
200 & 1.09266 & 0.04901\\
300 & 1.09870 & 0.04777\\
500 & 1.10544 & 0.04646\\ \hline
\(\infty\) & \(3.6/\pi\approx1.14592\) & \(1/(8\pi)\approx0.03979\)
\end{tabular}
\end{center}
The theorem gives no rate. The explicit bound
\(a_n\le1.14736\,g_n(1)\) of Theorem~\ref{thm:Ubound} is asymptotically
off by a factor \(1.0013\) only.
\end{remark}

\section{Conclusion and conjectures}\label{sec:open}

We have shown that the number of closed linear \(\lambda\)-terms in
\(\beta\)-normal form with \(n\) variable occurrences satisfies
\[
a_n\sim\frac{(n-1)!}{8\pi}\Bigl(\frac1{\ln2-\frac12}\Bigr)^n .
\]
Each factor has a structural origin. The growth rate is the inverse of the
instanton action \(A_0=\ln2-\frac12\), which is the radius of the Borel
transform of the linearised class (\(\lambda\)-chain terms). The factor
\(\frac1{2\pi}=\kappa_\infty\) comes from the turning point, through an exact
reduction to the logarithm of the Airy series. The factor
\(\frac14=\mu(0)\) is the transported WKB prefactor, which appears as an exact
conjugation \(T''=\mu\,\partial_s\,\mu^{-1}\) of the linear dynamics. The
proof rests on positivity in the variable \(v=1/r\), on an exact identity
whose remainder is either negative or of relative size \(O(m^{-2})\) with
respect to the reference profile, and on the dominance of the top
coefficient under the conjugated dynamics. The only computer-assisted input
is a finite bound on \(\sup_mU_m\) (Theorem~\ref{thm:Ubound},
Appendix~\ref{app:code}); it also gives the explicit bound
\(a_n\le1.14736\,g_n(1)\) for every \(n\ge2\).

We conclude with three conjectures.

\paragraph{Fixed \(w\).} All coefficientwise inequalities of
Sections~\ref{sec:pos}--\ref{sec:main} hold in the variable \(v\), and
Theorem~\ref{thm:main} only evaluates them at \(v=1\), that is \(w=0\).

\begin{conjecture}\label{conj:w}
For every fixed \(w\in[0,\frac14)\),
\[
\varphi_m(w)\sim\frac{\mu(w)}{2\pi}\,(m-1)!\,A(r)^{-m}
\qquad(m\to\infty),\qquad \mu(w)=\frac{r}{(1+r)^2},\quad r=\sqrt{1-4w},
\]
\end{conjecture}

\noindent Here \(A(r)=\frac{r^2}2-r+\ln(1+r)\) is the function of
Theorem~\ref{thm:action}, evaluated at \(r=\sqrt{1-4w}\); thus
\(A(r)=A_0-s(w)\). For \(w=0\) this is Theorem~\ref{thm:intro}. For \(w>0\) the proof of
Theorem~\ref{thm:main} should carry over by evaluating the inequalities at
\(v=1/r(w)\) and expanding around \(s(w)\) instead of \(0\): the singularity
at \(s=A_0\) is then at distance \(A_0-s(w)=A(r)\). We have not written the
details.

\paragraph{The \(n^{-1/3}\) correction.} Least-squares fits of
\(a_n/((n-1)!\lambda^n)\) on exact data up to \(n=260\), in the basis
\(n^{-j/3}\), are stable, whereas fits in integer powers of \(1/n\) do not
converge.

\begin{conjecture}\label{conj:sub}
\[
\frac{a_n}{(n-1)!\,\lambda^n}=\frac1{8\pi}\Bigl(1+b_1n^{-1/3}+O(n^{-2/3})\Bigr),
\qquad b_1=\kappa A_0^{1/3}\approx1.2599,\quad\kappa\approx2.1795 .
\]
\end{conjecture}

The fitted value \(\kappa=2.1795\pm0.0001\) is stable. The same extraction
for \(\varphi_n(w)\) at other values of \(w\) gives
\(b_1(w)/A(r)^{1/3}\approx2.179\), with \(r=\sqrt{1-4w}\), independent of \(w\). This suggests a
Stokes multiplier \(1+\kappa\,\varepsilon^{1/3}+\cdots\) generated at the
turning point by the corrections to the Airy equation. In the decomposition
\(a_n=\frac{a_n}{g_n(1)}\cdot g_n(1)\), the second factor contributes
\(c_f\approx1.570\) to \(b_1\) (Proposition~\ref{prop:class}), so the
convergence of \(a_n/g_n(1)\) would contribute about \(-0.31\). No closed
form for \(\kappa\) has been identified.

\paragraph{Depth.} View a term as a tree, with the root at depth \(0\), the
body of an abstraction and both children of an application one level
deeper. Let \(D(x,t)\) and \(E(x,t)\) be the generating functions, weighted
as \(F\), of the total depth of the variable occurrences and of all
\(3n-k-1\) nodes. Decomposing as in \eqref{eq:core} (the head of a spine
\(x\,N_1\cdots N_p\) lies at depth \(p\), and \(N_i\) at depth \(p-i+1\))
gives the linear equations
\[
D=\partial_t\bigl(D+x\partial_xF\bigr)+\frac{xt\,(F+D)}{(1-F)^2}+\frac{xt\cdot x\partial_xF}{(1-F)^3},
\qquad
E=\partial_t\bigl(E+\mathcal SF\bigr)+xt\Bigl[\frac{F+\mathcal SF}{(1-F)^3}+\frac{E}{(1-F)^2}\Bigr],
\]
with \(\mathcal S=3x\partial_x-t\partial_t-1\). They determine the exact mean
depths for closed terms; these values coincide with exhaustive enumeration
for \(n\le5\), and were computed up to \(n=140\) (at \(n=140\): \(168.41\)
for variables and \(138.10\) for nodes).

\begin{conjecture}\label{conj:depth}
In a uniformly random closed linear normal term with \(n\) variable
occurrences, the mean depth of a variable occurrence and the mean depth of
a node are asymptotically \(\alpha\,n\) and \(\beta\,n\), with
\(\alpha\approx1.18622\) and \(\beta\approx0.98086\).
\end{conjecture}

The constants come from extrapolating the exact values up to \(n=140\); two
fitting bases and two sets of points agree to about \(10^{-5}\). The
equations for \(D\) and \(E\) have the same principal part as
\eqref{eq:core}, so the methods of this paper should apply to them.

Finally, it would be desirable to bound \(\sup_mU_m\) without computer
assistance, and to determine the asymptotics of \(f_n(t)\) for fixed
\(t>0\).

\appendix
\section{Program for Theorem~\ref{thm:Ubound}}\label{app:code}

The following Python program computes the integer profiles
\(4^m\varphi_{m,\cdot}\) and \(4^mg_{m,\cdot}\), the exact maxima \(U_m\) by
integer cross-multiplication, checks \(a_1,\dots,a_7\) and the positivity of
all profiles, and checks every numerical inequality of the tail argument
and of the bootstrap in the proof of Theorem~\ref{thm:Ubound} in exact
rational arithmetic (Python \texttt{Fraction}); no floating-point number
enters the proof, floats being used only to display \(U_m\).
Polynomial products use \texttt{python-flint} when it is installed (about 8
minutes for \(M=700\)) and an exact pure-Python Kronecker substitution
otherwise (about an hour). The output below was produced with
\texttt{python-flint}; run with \(M=300\), the pure-Python backend prints the
same values of \(U_{100},U_{200},U_{300}\), the same location of the
maxima and the same \(\max_{m\le300}U_m\). The mathematical content of the program
(the recurrences of Lemma~\ref{lem:pos} and Proposition~\ref{prop:class})
can be checked independently of its execution.

{\footnotesize
\begin{verbatim}
#!/usr/bin/env python3
"""Exact verification for Theorem 7.2 (bound on the depth ratio U_m).

All profile quantities are integers: P_m = 4^m phi_{m,.}, F_m = 4^m g_{m,.}.
The tail constants are checked in exact rational arithmetic (no floats).
Polynomial products use python-flint if available, otherwise an exact
pure-Python Kronecker substitution (slower). Runtime for M=700: about
8 minutes with python-flint, about an hour without.
"""
import sys, math
from fractions import Fraction as Fr
M = int(sys.argv[1]) if len(sys.argv) > 1 else 700

try:
    import flint
    def mul(p, q):                       # exact product of integer polynomials
        return [int(c) for c in (flint.fmpz_poly(p) * flint.fmpz_poly(q)).coeffs()]
    BACKEND = "python-flint"
except ImportError:
    def mul(p, q):                       # Kronecker substitution, nonnegative integers
        bits = max(max(p).bit_length(), 1) + max(max(q).bit_length(), 1) \
               + max(len(p), len(q)).bit_length() + 1
        pack = lambda v: sum(c << (bits * i) for i, c in enumerate(v))
        x, mask, out = pack(p) * pack(q), (1 << bits) - 1, []
        for _ in range(len(p) + len(q) - 1):
            out.append(x & mask); x >>= bits
        return out
    BACKEND = "pure Python"
print("backend:", BACKEND)

def add_into(acc, v, factor=1):
    for i, c in enumerate(v):
        acc[i] += factor * c

# Lemma 5.1 (positive recurrence), scaled by 4^m:
#   X_m[d] = 4(3m-4-d)(P_{m-1}[d] + X_{m-1}[d-1]),  P_m = X_m + Q_m,
#   Q_m = sum_{i=1}^{m-1} P_i * P_{m-i}  (depth convolution).
P = {1: [2, 2]}; X = [2, 2]
F = [20, 24, 4]                      # 16 * g_2 = 16 * phi_2 = 16*(5/4, 3/2, 1/4)
U = {}; argmax = {}; a = {1: 1}
for m in range(2, M + 1):
    Q = [0] * (m + 1)
    for i in range(1, (m - 1) // 2 + 1):
        add_into(Q, mul(P[i], P[m - i]), 2)
    if m % 2 == 0:
        add_into(Q, mul(P[m // 2], P[m // 2]))
    Pm1 = P[m - 1] + [0, 0]; Xm1 = X + [0, 0]
    X = [4 * (3*m - 4 - d) * (Pm1[d] + (Xm1[d-1] if d >= 1 else 0)) for d in range(m + 1)]
    P[m] = [X[d] + Q[d] for d in range(m + 1)]
    if m >= 3:  # g_m = T'' g_{m-1}:  F[d] = 4((3m-3-d)F_{m-1}[d] + (3m-4-d)F_{m-1}[d-1])
        Fm1 = F + [0]
        F = [4 * ((3*m - 3 - d) * Fm1[d] + ((3*m - 4 - d) * Fm1[d-1] if d >= 1 else 0))
             for d in range(m + 1)]
    assert all(x > 0 for x in F) and all(x >= 0 for x in P[m])
    b = 0
    for d in range(1, m + 1):              # exact comparison P[d]/F[d] > P[b]/F[b]
        if P[m][d] * F[b] > P[m][b] * F[d]:
            b = d
    argmax[m] = b; U[m] = Fr(P[m][b], F[b])
    s = sum(P[m]); assert s % 4**m == 0; a[m] = s // 4**m
    if m % 100 == 0:
        print(f"m={m}: U_m={float(U[m]):.10f}, argmax d={b}", flush=True)

assert [a[n] for n in range(1, 8)] == [1, 3, 26, 367, 7142, 176766, 5304356]
Ubar = max(U.values()); mbar = max(U, key=lambda m: U[m])
print("argmax d = 0 for 5<=m<=M:", all(argmax[m] == 0 for m in range(5, M + 1)))
print(f"max_(m<=M) U_m = U_{mbar} = {float(Ubar):.12f}")

# Tail bound for m > 700 (proof of Theorem 7.2), in exact rational arithmetic:
# every inequality used in the proof is checked with Fractions.
assert M >= 700, "the bootstrap of Theorem 7.2 starts from m = 700"
gam = lambda k: Fr(3*k - 3, 2*k - 4)                       # gamma_k
D = Fr                                                     # Fr("0.37568") is exact
rho, T2c, C1, q = D("0.37568"), D("0.31340"), D("1.0040"), D("0.8914")
R0, eta, Ustar = D("0.70"), D("0.001437"), D("1.14736")
checks = [
  ("gamma_561/4 <= 0.37568",                gam(561) / 4 <= rho),
  ("(5/36) gamma_700 gamma_701 <= 0.31340", Fr(5, 36) * gam(700) * gam(701) <= T2c),
  ("2*0.31340/(1-0.37568) <= 1.0040",       2 * T2c / (1 - rho) <= C1),
  ("gamma_351 <= 1.51",                     gam(351) <= D("1.51")),
  ("5 - 10/701 >= 4.98",                    5 - Fr(10, 701) >= D("4.98")),
  ("(5/36)*4.98*701/700 <= 0.70",           Fr(5, 36) * D("4.98") * Fr(701, 700) <= R0),
  ("(1.51^(5/2)/4.98)^(1/5) <= 0.8914",     D("1.51")**5 <= (D("4.98") * q**5)**2),
  ("0.70*0.8914^701/(1-0.8914) <= 1e-34",   R0 * q**701 / (1 - q) <= Fr(1, 10**34)),
  ("1.0040/699 + 1e-34 <= 0.001437",        C1 / 699 + Fr(1, 10**34) <= eta),
  ("Ubar <= 1.1454597",                     Ubar <= D("1.1454597")),
  ("Ubar + 1.14736^2 * 0.001437 <= 1.14736", Ubar + Ustar**2 * eta <= Ustar),
]
for name, ok in checks:
    print(f"{name}: {ok}")
assert all(ok for _, ok in checks)
\end{verbatim}
}

\noindent Output (\texttt{python3 verify\_U.py 700}):
{\footnotesize
\begin{verbatim}
backend: python-flint
m=100: U_m=1.1426714638, argmax d=0
m=200: U_m=1.1443090123, argmax d=0
m=300: U_m=1.1448479108, argmax d=0
m=400: U_m=1.1451160860, argmax d=0
m=500: U_m=1.1452765871, argmax d=0
m=600: U_m=1.1453834202, argmax d=0
m=700: U_m=1.1454596478, argmax d=0
argmax d = 0 for 5<=m<=M: True
max_(m<=M) U_m = U_700 = 1.145459647843
gamma_561/4 <= 0.37568: True
(5/36) gamma_700 gamma_701 <= 0.31340: True
2*0.31340/(1-0.37568) <= 1.0040: True
gamma_351 <= 1.51: True
5 - 10/701 >= 4.98: True
(5/36)*4.98*701/700 <= 0.70: True
(1.51^(5/2)/4.98)^(1/5) <= 0.8914: True
0.70*0.8914^701/(1-0.8914) <= 1e-34: True
1.0040/699 + 1e-34 <= 0.001437: True
Ubar <= 1.1454597: True
Ubar + 1.14736^2 * 0.001437 <= 1.14736: True
\end{verbatim}
}

\end{document}